\documentclass[12pt,a4paper]{amsart}

\usepackage{amssymb}
\usepackage{amsfonts}
\usepackage{amsthm}
\usepackage{amsmath} 
\usepackage{enumerate} 
\usepackage{tikz} 
\usepackage{mathtools} 
\usepackage[latin1]{inputenc}
\usepackage[british]{babel}
\usepackage[colorlinks]{hyperref}	

\hypersetup{
    colorlinks,
    citecolor=green,
    linkcolor=blue
}	

\numberwithin{equation}{section}

\theoremstyle{plain}
	\newtheorem{Thm}{Theorem}[section]
	\newtheorem{Lem}[Thm]{Lemma}
	\newtheorem{Conj}[Thm]{Conjecture}

\theoremstyle{definition}
	\newtheorem{Def}[Thm]{Definition}
	\newtheorem{Ex}{Example}
	
\theoremstyle{remark}

	\newtheorem{Case}{Case}
	\newtheorem{Step}{Step}

\newcommand{\edges}{E} 
\newcommand{\vertices}{V} 
\DeclareMathOperator\ldeg{ldeg}
\DeclareMathOperator\gdeg{gdeg}
\DeclareMathOperator\mdeg{mdeg}
\newcommand{\nats}{\mathbb{N}}
\newcommand{\ints}{\mathbb{Z}}
\DeclareMathOperator\res{res}

\newcommand{\mon}[1]{\mathbf{x}^{\mathbf{#1}}}

\newcommand{\bei}{\mathcal{J}_G} 
\newcommand{\pbei}{\mathcal{I}_G} 
\newcommand{\gb}[2]{\mathcal{G}_{#1,#2}} 
\newcommand{\gbp}[1]{\mathcal{G}_{\succ}({#1})} 

\newcommand{\ugb}[1]{\mathcal{U}(#1)} 
\newcommand{\gra}[1]{\mathrm{Gr}(#1)} 
\newcommand{\sset}[1]{\mathcal{S}(#1)} 

\title{The universal Gr\"obner basis of a binomial edge ideal}

\author{Mourtadha Badiane}
\address{Mourtadha Badiane, National University of Ireland, Galway}
\email{m.badiane1@nuigalway.ie}

\author{Isaac Burke}
\thanks{The second author was supported by the Irish Research Council and the Hardiman Scholarship Scheme}
\address{Isaac Burke, National University of Ireland, Galway}
\email{i.burke1@nuigalway.ie}

\author{Emil Sk\"oldberg}
\address{Emil Sk\"oldberg, National University of Ireland, Galway}
\email{emil.skoldberg@nuigalway.ie}

\date{\today}



\begin{document}

\maketitle

\begin{abstract}
  We show that the universal Gr\"obner basis and the Graver basis of a binomial edge ideal coincide. We provide a description for this basis set in terms of certain paths in the underlying graph. We conjecture a similar result for a parity binomial edge ideal and prove this conjecture for the case when the underlying graph is the complete graph.
\end{abstract}
%
%
%
%
\section{Introduction}
%
%
%
%
For $n \in \nats_{>0}$, $[n]:=\{1,...,n\}$. Let $G$ be a simple graph on the vertex set $[n]$, that is, $G$ has no loops and no multiple edges. Let $\edges(G)$ denote the edge set of $G$. Let $F$ be a field and let $S=F[x_1,\ldots,x_n,y_1,\ldots,y_n]$ be the polynomial ring in $2n$ variables. The binomial edge ideal of $G$ was introduced and studied independently by Herzog, Hibi, Hreinsd\'ottir, Kahle and Rauh~\cite{hhhkr:binomial_edge_ideals} and
Ohtani~\cite{ohtani:binomial_edge_ideals}.
\begin{Def}
The \emph{binomial edge ideal} of $G$ is
\begin{equation}
		\bei:=
		\langle x_iy_j-x_jy_i : 
		\ \{i,j\} \in \edges(G) 
		\rangle \subseteq S.
\end{equation}
\end{Def}
The parity binomial edge ideal of $G$ was introduced and studied by Kahle, Sarmiento and Windisch~\cite{ksw:parity_binomial_edge_ideals} but had previously been examined by Herzog, Macchia, Madani and Welker~\cite{hmmw:ideal_of_orthogonal_representations}. 
\begin{Def}
The \emph{parity binomial edge ideal} of $G$ is
\begin{equation}\label{eq:ideal}
		\pbei:=
		\langle x_ix_j-y_iy_j : 
		\{i,j\} \in \edges(G) 
		\rangle \subseteq S.
\end{equation}
\end{Def}
These ideals appear in various settings and applications in mathematics and statistics and belong to an important class of binomial ideals which may be defined as follows. If we let $R=F[x_1,\ldots,x_n]$ then an ideal $I$ of $R$ is a \emph{pure difference ideal} (also known in the literature as a \emph{pure binomial ideal}) if $I$ is generated by differences of monic monomials i.e.\ binomials of the form $\mon{u}-\mon{v}$ with $u,v \in \nats^n$. There are several well-known distinguished subsets of binomials in such an ideal $I$, two of which we now mention. A binomial $\mon{u}-\mon{v} \in I$ is called \emph{primitive} if there exists no other binomial $\mon{u'}-\mon{v'} \in I$ such that $\mon{u'}$ divides $\mon{u}$ and $\mon{v'}$ divides $\mon{v}$. The set of primitive binomials in $I$ is called the \emph{Graver basis} of $I$ and denoted $\gra{I}$. The union of all of the reduced Gr\"obner bases of $I$ is called the \emph{universal Gr\"obner basis} of $I$ and denoted $\ugb{I}$. Graver bases were originally defined for toric ideals by Sturmfels~\cite{sturmfels:groebner_bases_convex_polytopes}. Charalambous, Thoma and Vladoiu~\cite{ctv:binomial_fibers_indispensable_binomials} recently generalised the concept to an arbitrary pure difference ideal $I$, showing in particular that $\gra{I}$ is finite and includes $\ugb{I}$ as a subset. 

One open problem that arises in the literature is providing a combinatorial characterisation of toric ideals for which the universal Gr\"obner basis and the Graver basis are equal (many examples have been discovered, see Petrovi{\'c} et al.~\cite{ptv:bouquet_algebra_toric_ideals} and references therein). We consider this problem for certain classes of pure difference ideals which are not lattice ideals. In particular, we show that $\ugb{\bei}=\gra{\bei}$ and provide a description for this basis set in terms of certain paths in $G$. We conjecture a similar result for $\pbei$ and prove this conjecture for the case when $G$ is the complete graph.

\subsection{Preliminaries}

Throughout the paper we assume that $G$ is finite, undirected and connected. For any $W \subseteq [n]$, let $G[W]$ denote the induced subgraph on $W$, and for a sequence of vertices $\pi=(i_0,...,i_r) \in [n]^{r+1}$, $G[\pi]:=G[\{i_0,...,i_r\}]$. A $(v,w)${\it-path} of {\it length} $r$ is a sequence of vertices $v=i_0,i_1,\ldots,i_r=w$ such that $\{i_k,i_{k+1}\} \in \edges(G)$ for all $k=0,\ldots,r-1$. The path is {\it odd (even)} if its length is odd (even). The {\it interior} of a $(v,w)$-path $\pi=(i_0,\ldots,i_r)$ is the set $\mathrm{int}(\pi)=\{i_0,\ldots,i_r\} \setminus \{v,w\}$. The {\it inverse} $\pi^{-1}$ of a $(v,w)$-path $\pi=(i_0,\ldots,i_r)$ is the $(w,v)$-path $(i_r,i_{r-1},\ldots,i_0)$. For the vertex set of a graph $H$ we sometimes use the notation $\vertices(H)$. For a monomial $\mon{u}=x_1^{d_1}y_1^{e_1} \cdots x_n^{d_n}y_n^{e_n}$ in $S$ the set $\{i:d_i\neq0 \ \mathrm{or} \ e_i\neq0\} \subseteq [n]$ is denoted by $\vertices(\mon{u})$.
%
%
%
%
\section{Binomial Edge Ideals}
%
%
%
%
In this section we will use two different gradings on $S$, the first is the $\nats^2$-grading by considering the letter of a variable, so we let $\ldeg(x_i) = (1,0)$ and $\ldeg(y_i) = (0,1)$ for all $i \in [n]$. The second is the $\nats^n$-grading which considers the vertex of a variable and we set $\gdeg(x_i) = \gdeg(y_i) = \mathbf{e}_i$ for all $i \in [n]$, where $\mathbf{e}_i$ is the $i$th standard basis vector in $\nats^n$. The ideal $\bei$ is homogeneous with respect to both of these gradings and we combine them into what we call the \emph{multidegree} of a monomial
$\mdeg(\mon{u}):=(\ldeg(\mon{u}),\gdeg(\mon{u})) \in {\nats}^2 \times {\nats}^{n}$. 

We now recall the definition of admissible paths and the description of the Gr\"obner basis of $\bei$ with respect to the lexicographic order which was independently obtained by Herzog et al.~\cite{hhhkr:binomial_edge_ideals} and Ohtani~\cite{ohtani:binomial_edge_ideals}.

\begin{Def}
Fix a permutation $\sigma \in S_n$ of $[n]$ and let $i,j \in [n]$ satisfy $\sigma^{-1}(i) < \sigma^{-1}(j)$. An $(i,j)$-path $\pi=(i_0,\ldots,i_r)$ in $G$ is called \emph{$\sigma$-admissible}, if
\begin{enumerate}[(i)]
\item $i_k \neq i_l$ if $k \neq l$;
\item $j_0, \ldots, j_s$ is not a path from $i$ to $j$ for any proper
  		  subset $\{j_0, \ldots, j_s\}$ of $\{i_0, \ldots, i_r\}$;
\item for each $k = 1, \ldots, r-1$, either $\sigma^{-1}(i_k) < \sigma^{-1}(i)$ or
      $\sigma^{-1}(i_k) > \sigma^{-1}(j)$.
\end{enumerate}
\end{Def}

Given a $\sigma$-admissible $(i,j)$-path $\pi=(i_0,\ldots,i_r)$ in $G$, where $\sigma^{-1}(i) < \sigma^{-1}(j)$,
\begin{equation*}
		u_{\pi} := 
		\prod_{\sigma^{-1}(i_k) < \sigma^{-1}(i)} y_{i_k}
		\prod_{\sigma^{-1}(i_k) > \sigma^{-1}(j)} x_{i_k}. 
\end{equation*}

\begin{Thm}[\cite{hhhkr:binomial_edge_ideals}, \cite{ohtani:binomial_edge_ideals}]\label{thm:gb}
The set of binomials
\begin{equation*}
		\gb{G}{\sigma} := \ 
		\bigcup_{\mathclap{\sigma^{-1}(i) < \sigma^{-1}(j)}} \ \{ u_{\pi} (x_iy_j - x_jy_i) : 
		\text{$\pi$ is a $\sigma$-admissible $(i,j)$-path in $G$} \}
\end{equation*}
is the reduced Gr\"obner basis of $\bei$ w.r.t.\ the lexicographic monomial order on $S$ induced by $x_{\sigma(1)} \succ \dots \succ x_{\sigma(n)} \succ y_{\sigma(1)} \succ \dots \succ y_{\sigma(n)}$.
\end{Thm}

Our first result is a characterisation of the binomials in $\bei$. For this we need to introduce the following notations. We denote by $d_G(v,w)$ the length of a shortest $(v,w)$-path in $G$. For a monomial $\mon{u}=x_1^{d_1}y_1^{e_1} \cdots x_n^{d_n}y_n^{e_n} \in S$, we sometimes use the notation $\deg_{x_i}(\mon{u})$ for $d_i$ and $\deg_{y_i}(\mon{u})$ for $e_i$. For an induced subgraph $H$ of $G$, we define the \emph{restriction} of $\mon{u}$ to $H$ to be $\res_{H}(\mon{u}) = \Pi_{i \in \vertices(H)}x_i^{d_i}y_i^{e_i}$. 

\begin{Lem}\label{lemma:multi-homogeneous}
Let $\mon{u} - \mon{v}$ be a multi-homogeneous binomial such that $G[\vertices(\mon{u})]$ is a connected graph. Then $\mon{u} - \mon{v} \in \bei$.
\proof
Suppose that $m = \mon{u}$ is a monomial such that there is a pair of indices $i < j$ with $\deg_{x_i}(m) \geq 1$ and $\deg_{y_j}(m) \geq 1$. Let $G' =  G[\vertices(m)]$. We can assume that $i$ and $j$ are chosen such that $d_{G'}(i,j)$ is minimal. Now we consider an $(i,j)$-path $\pi=(i_0,\ldots,i_r)$ of minimal length in $G'$. Since $\pi$ is of minimal length we can conclude that $i_k \neq i_l$ for $k \neq l$ and that no proper subset $\{j_0, \ldots, j_s\}$ of $\{i_0, \ldots, i_r\}$ is a path from $i$ to $j$. Suppose that there is a $k$ such that $i < i_k < j$, then either $\deg_{x_{i_k}}(m) \geq 1$, in which case $d_{G'}(i_k, j) < d_{G'}(i,j)$ which contradicts the minimality of $d_{G'}(i,j)$, or $\deg_{y_{i_k}}(m) \geq 1$, in which case $d_{G'}(i, i_k) < d_{G'}(i,j)$ which again contradicts the minimality of $d_{G'}(i,j)$. We may thus conclude that $\pi$ is a $\sigma$-admissible $(i,j)$-path in $G'$, where $\sigma=\text{id}$, the identity permutation in $S_n$.

Now we consider the vertex $i_k$ on $\pi$. By the minimality of $d_{G'}(i,j)$, if $i_k < i$ then $\deg_{x_{i_k}}(m) = 0$ and if $i_k > j$ then $\deg_{y_{i_k}}(m) = 0$. We may thus conclude that $u_{\pi} x_iy_j$ divides $m$, and therefore $m$ is reducible with respect to $\gb{G}{\mathrm{id}}$. This shows that an irreducible monomial of the same multidegree as $\mon{u}$ has the form $\mon{w} = y_{i_1}^{e_1}\cdots y_{i_k}^{e_k}x_{i_k}^{d_k}\cdots x_{i_l}^{d_l}$ where $i_1 < i_2 < \cdots <i_l$. Since there is only one such monomial in a given multidegree, we can conclude that $\mon{u} - \mon{v}$ reduces to zero with respect to $\gb{G}{\mathrm{id}}$ and thus $\mon{u} - \mon{v} \in \bei$.
\qed
\end{Lem}

\begin{Lem}\label{lemma:membership}
A multi-homogeneous binomial $\mon{u} - \mon{v}$ lies in $\bei$ if and only if $\mdeg(\res_{C}(\mon{u})) = \mdeg(\res_{C}(\mon{v}))$ for all connected components $C$ of $G[\vertices(\mon{u})]$.
\proof 
Let $\mon{u} - \mon{v} \in \bei$, then we can write
\begin{equation*}
		\mon{u} - \mon{v} = \sum_{k} 
		\mathbf{x}^{\mathbf{w}_k} 
		(x_{i_k}y_{j_k}-x_{j_k}y_{i_k})
\end{equation*}
where $\{i_k, j_k\} \in \edges(G)$ for all $k$. Let $C$ be a component in $G[\vertices(\mon{u})]$, then by restricting to $C$ we get
\begin{equation*}
		\res_{C}(\mon{u}) - \res_{C}(\mon{v}) = 
		\sum_{k, i_k \in \vertices(C)} \res_{C} 
		(\mathbf{x}^{\mathbf{w}_k})
		(x_{i_k}y_{j_k}-x_{j_k}y_{i_k})
\end{equation*}
since $j_k \in \vertices(C)$ if and only if $i_k \in \vertices(C)$. Thus we see that $\res_C(\mon{u}) -  \res_C(\mon{v}) \in \bei$, and therefore $\mdeg(\res_{C}(\mon{u})) = \mdeg(\res_{C}(\mon{v}))$.

For the converse, suppose $\mon{u} - \mon{v}$ satisfies that $\mdeg(\res_{C}(\mon{u})) = \mdeg(\res_{C}(\mon{v}))$ for all connected components $C_1, \ldots C_r$ of $G[\vertices(\mon{u})]$, then, by Lemma~\ref{lemma:multi-homogeneous}, $\res_{C_i}(\mon{u}) - \res_{C_i}(\mon{v}) \in \bei$ for all $r$, and we can write
\begin{equation*}
 	 	\mon{u} - \mon{v} = \sum_{i=1}^{r}
 		\frac{\mon{u} \prod_{j=1}^{i-1}\res_{C_i}(\mon{v})}
  		{\prod_{j=1}^{i}\res_{C_i}(\mon{u})}
  		(\res_{C_i}(\mon{u}) - \res_{C_i}(\mon{v})),
\end{equation*}
which shows that $\mon{u} - \mon{v} \in \bei$.
\qed
\end{Lem}

\begin{Def}
An $(i,j)$-path $\pi=(i_0,\ldots,i_r)$ in $G$ is called \emph{weakly admissible} if it satisfies conditions (i) and (ii) of the definition of a $\sigma$-admissible path.
\end{Def}

Given a weakly admissible $(i,j)$-path $\pi=(i_0,\ldots,i_r)$ in $G$,
\begin{equation*}
		\mathcal{S}_{\pi} := 
		\{ t_{i_1}t_{i_2}\cdots t_{i_{r-1}}(x_iy_j - x_jy_i) : 
		t_k \in \{x_k, y_k\} \}
\end{equation*}
and $\sset{\bei} := \bigcup_{\pi} \mathcal{S}_{\pi} \setminus \{0\}$ where $\pi$ runs over all weakly admissible paths in $G$. Notice that if $\pi$ is an $(i,i)$-path in $G$, then $\pi$ is weakly admissible if and only if $\pi$ is the path $(i)$ of length 0, in which case $\mathcal{S}_{\pi}=\{0\}$.

\begin{figure}[t]
	\begin{tikzpicture}[xscale=0.5,yscale=0.5]
		\node [label={[label distance=1pt]180:$3$},fill, circle, inner sep=2pt](a) at (0,0) {};
		\node [label={[label distance=1pt]180:$1$},fill, circle, inner sep=2pt](b) at (0,2) {};
		\node [label={[label distance=1pt]90:$2$},fill, circle, inner sep=2pt](c) at (2,1) {};
		\node [label={[label distance=1pt]90:$4$},fill, circle, inner sep=2pt](d) at (4,1) {};
		\draw(a)--(b)--(c)--(a);
		\draw(c)--(d);
	\end{tikzpicture}
	\caption{See Examples~\ref{ex:bei} and \ref{ex:pbei}.} \label{fig:graph}
\end{figure}

\begin{Ex}\label{ex:bei}
Let $G$ be the graph in Figure~\ref{fig:graph}. The weakly admissible paths in $G$ are the paths $(1), \ (2), \ (3), \ (4), \ (1,2), \ (1,3), \ (2,3), \ (2,4)$, $(1,2,4)$ and $(3,2,4)$, together with their inverses. Hence $|\sset{\bei}|=16$.
\end{Ex}

\begin{Thm}
The sets $\sset{\bei}$, $\ugb{\bei}$ and $\gra{\bei}$ coincide.
\proof
We prove the theorem in three steps; the containments $\sset{\bei} \subseteq \ugb{\bei}$, $\ugb{\bei} \subseteq \gra{\bei}$ and $\gra{\bei} \subseteq \sset{\bei}$.

\begin{Step}
$\sset{\bei} \subseteq \ugb{\bei}$: Let $\pi=(i_0,\ldots,i_r)$ be a weakly admissible $(i,j)$-path in $G$ and let $f = t_{i_1}t_{i_2}\cdots t_{i_{r-1}}(x_iy_j - x_jy_i)$ be a corresponding binomial in $\sset{\bei}$. Now let $\sigma \in S_n$ be a permutation such that $\sigma^{-1}(i) < \sigma^{-1}(j)$, $\sigma^{-1}(i_k) < \sigma^{-1}(i)$ for all $k$ such that $t_{i_k} = y_{i_k}$ and $\sigma^{-1}(i_k) > \sigma^{-1}(j)$ for all $k$ such that $t_{i_k} = x_{i_k}$. Then $f \in \gb{G}{\sigma}$ and thus $f \in \ugb{\bei}$.
\end{Step}

\begin{Step}
$\ugb{\bei} \subseteq \gra{\bei}$: \cite[Proposition 4.2]{ctv:binomial_fibers_indispensable_binomials}.
\end{Step}

\begin{Step}
$\gra{\bei} \subseteq \sset{\bei}$: Let $\mon{u} - \mon{v}$ be a primitive binomial in $\bei$ and  let $C$ be a component of $G[\vertices(\mon{u})]$. By Lemma~\ref{lemma:membership} we have that $\res_{C}(\mon{u}) - \res_{C}(\mon{v}) \in \bei$, so by primitivity we can conclude that $C := G[\vertices(\mon{u})]$ is connected.

Since $\mathbf{u} \neq \mathbf{v}$ we can choose $i,j \in \vertices(\mon{u})$ such that $d_C(i,j)$ is minimal among all pairs $i, j$ with $\deg_{x_i}(\mon{u}) > \deg_{x_i}(\mon{v})$ and $\deg_{y_j}(\mon{u}) > \deg_{y_j}(\mon{v})$. Now let $\pi=(i_0,\ldots,i_r)$ be an $(i,j)$-path in $C$ of minimal length. Suppose that there is a $k \in \{1,\ldots,r-1\}$ such that $\deg_{x_{i_k}}(\mon{u}) \neq \deg_{x_{i_k}}(\mon{v})$. Then either $\deg_{x_{i_k}}(\mon{u}) > \deg_{x_{i_k}}(\mon{v})$, in which case $d_C(i_k,j)$ would contradict the minimality of $d_C(i,j)$, or $\deg_{x_{i_k}}(\mon{u}) < \deg_{x_{i_k}}(\mon{v})$, in which case by homogeneity we would have $\deg_{y_{i_k}}(\mon{u}) > \deg_{y_{i_k}}(\mon{v})$ and thus $d_C(i,i_k)$ would contradict the minimality of $d_C(i,j)$. So for $k \in \{1,\ldots,r-1\}$ we have $\deg_{x_{i_k}}(\mon{u}) = \deg_{x_{i_k}}(\mon{v})$ and hence by homogeneity $\deg_{y_{i_k}}(\mon{u}) = \deg_{y_{i_k}}(\mon{v})$. We can then, for $k = 1,\ldots,r-1$, let $z_{i_k} \in \{x_{i_k}, y_{i_k}\}$ such that $z_{i_k}$ divides $\mon{u}$ and thus also $\mon{v}$. Then $z_{i_1}\cdots z_{i_{r-1}} (x_iy_j-x_jy_i) \in \bei$ with $z_{i_1}\cdots z_{i_{r-1}} x_iy_j | \mon{u}$ and $z_{i_1}\cdots z_{i_{r-1}} x_jy_i | \mon{v}$, which implies that $\mon{u} - \mon{v} = z_{i_1}\cdots z_{i_{r-1}} (x_iy_j-x_jy_i) \in \sset{\bei}$.
\qed
\end{Step}
\end{Thm}
%
%
%
%
\section{Parity Binomial Edge Ideals}
%
%
%
%
In this section we will use two different gradings on $S$, but not exactly as in the previous section. The first grading is the $\ints^2_2$-grading by considering the letter of a variable, so we let $\ldeg(x_i)=(1,0) \in \ints^2_2$ and $\ldeg(y_i)=(0,1) \in \ints^2_2$ for all $i \in [n]$. The second is the $\nats^n$-grading as in the previous section. The ideal $\pbei$ is homogeneous with respect to both of these gradings and we combine them into what we call the {\it multidegree} of a monomial $\mdeg(\mon{u}):=(\ldeg(\mon{u}),\gdeg(\mon{u})) \in \ints^2_2 \times \nats^n$.

\begin{Lem}[\cite{ksw:parity_binomial_edge_ideals}]\label{lemma:ksw15}
Let $\pi=(i_0,\ldots,i_r)$ be an $(i,j)$-path in $G$ and $t_k \in \{x_k,y_k\}$ arbitrary. If $\pi$ is odd, then
\begin{equation*}
		(x_ix_j-y_iy_j)\prod_{k \in \mathrm{int}(\pi)}
		\limits t_k \in \mathcal{I}_G.
\end{equation*}
If $\pi$ is even, then
\begin{equation*}
		(x_iy_j-y_ix_j)\prod_{k \in \mathrm{int}(\pi)}
		\limits t_k \in \mathcal{I}_G.
\end{equation*}
\end{Lem}

\begin{Def}
An $(i,j)$-path $\pi$ in $G$ is called \emph{minimal}, if
\begin{enumerate}[(i)]
\item for no $k \in \mathrm{int}(\pi)$ there is an $(i,j)$-path with the same parity as $\pi$ in $G[\pi \setminus \{k\}]$;
\item there is no shorter $(i,j)$-path $\pi'$ in $G$ satisfying $\mathrm{parity}(\pi')=\mathrm{parity}(\pi)$ and $\mathrm{int}(\pi')=\mathrm{int}(\pi)$.
\end{enumerate} 
\end{Def}

For a minimal $(i,j)$-path $\pi$ in $G$, we define a set of binomials $\mathcal{S}_\pi$ as follows. If $\pi$ is odd, then $\mathcal{S}_\pi:=\mathcal{S}^{+}_{\pi,o} \bigcup \mathcal{S}^{-}_{\pi,o}$ where
\begin{align*}
		\mathcal{S}^{+}_{\pi,o} 
		&:=\{(x_ix_j-y_iy_j)\prod_{k \in \mathrm{int}(\pi)}
		\limits t_k:t_k \in \{x_k,y_k\}\}, \\
		\mathcal{S}^{-}_{\pi,o} 
		&:=\{(y_iy_j-x_ix_j)\prod_{k \in \mathrm{int}(\pi)}
		\limits t_k:t_k \in \{x_k,y_k\}\}.
\end{align*}
If $\pi$ is even, then $\mathcal{S}_\pi:=\mathcal{S}_{\pi,e}$ where
\begin{equation*}
		\mathcal{S}_{\pi,e}:=
		\{(x_iy_j-y_ix_j)\prod_{k \in \mathrm{int}(\pi)}
		\limits t_k:t_k \in \{x_k,y_k\}\}.
\end{equation*}
$\sset{\pbei}:=\bigcup_\pi \mathcal{S}_\pi \setminus \{0\}$ where $\pi$ runs over all minimal paths in $G$. Notice that if $\pi$ is an even $(i,i)$-path in $G$, then $\pi$ is minimal if and only if $\pi$ is the path $(i)$ of length 0, in which case $\mathcal{S}_\pi=\{0\}$.

\begin{Ex}\label{ex:pbei}
Let $G$ be the graph in Figure~\ref{fig:graph}. The minimal paths in $G$ are the paths
\begin{align*}
&(1), \ (2), \ (3), \ (4),\\
&(1,2), \ (1,3), \ (2,3), \ (2,4),\\
&(1,2,3), \ (1,2,4), \ (1,3,2), \ (2,1,3), \ (3,2,4),\\
&(1,2,3,1), \ (1,3,2,4), \ (2,1,3,2), \ (3,1,2,3), \ (3,1,2,4),\\
&(2,1,3,2,4), \ (2,3,1,2,4)\\
&\text{and } (4,2,1,3,2,4),
\end{align*}
together with their inverses. Hence $|\sset{\pbei}|=92$.
\end{Ex}

Given a graph $G$, it is clear that the set of its weakly admissible paths is a subset of the set of its minimal paths. If $\pi$ is a minimal $(i,j)$-path in $G$ which is not a weakly admissible path, $\pi$ contains repeated vertices or $G[\pi]$ contains an $(i,j)$-path $\pi'$ of parity opposite to that of $\pi$ such that $\mathrm{int}(\pi') \subsetneq \mathrm{int}(\pi)$. 
 
\begin{Conj}\label{conjecture:one}
The sets $\sset{\pbei}$, $\ugb{\pbei}$ and $\gra{\pbei}$ coincide.
\end{Conj}

We have partially tested Conjecture~\ref{conjecture:one} for small graphs using the software \textsc{gfan} \cite{gfan} and have found no counterexamples so far. It must be said, however, that we are not currently aware of any algorithm for computing the Graver basis of an arbitrary pure difference ideal. The main result of this section is a proof that Conjecture~\ref{conjecture:one} holds when $G$ is the complete graph $K_n$ on the vertex set $[n]$. The rest of the section is arranged as follows. In Lemmas~\ref{lemma:minimal} through \ref{lemma:basis} we describe a reduced Gr\"obner basis of $\mathcal{I}_{K_n}$. In Lemmas~\ref{lemma:normalform0} through \ref{lemma:normalform4} we characterise the binomials in $\mathcal{I}_{K_n}$. In Theorem~\ref{theorem:parity} the main result is proved.

\begin{Lem}\label{lemma:minimal}
The minimal paths in $K_n$ are all
\begin{equation*}
(i), \ (i,j), \ (i,k,j) \ \text{and } (i,k,l,i)
\end{equation*}
where $i,j,k$ and $l$ are distinct elements of $[n]$.
\proof
Let $\pi=(i_0,\ldots,i_r)$ be an $(i,j)$-path in $K_n$. Suppose that $\pi$ is odd. If $\mathrm{int}(\pi)=\varnothing$ then $\pi$ is necessarily of the form $(i,j,i,j,\ldots,i,j)$ and thus minimal if and only if $\pi=(i,j)$. If $\mathrm{int}(\pi) \neq \varnothing$ then there are two cases: $i=j$ and $i \neq j$. If $i=j$ then $|\mathrm{int}(\pi)| \geq 2$ so let $i_{s_1} \neq i_{s_2} \in \mathrm{int}(\pi)$ and notice that $K_n$ contains the odd path $\pi'=(i_0,i_{s_1},i_{s_2},i_0$). Now $K_n[\pi'\setminus \{i_{s_1}\}] \cong K_n[\pi'\setminus \{i_{s_2}\}] \cong K_2$ which does not contain an odd cycle, hence $\pi'$ is minimal. It follows that $\pi$ is minimal if and only if $\pi=(i,k,l,i)$, where $\{k,l\}=\mathrm{int}(\pi)$. If $i \neq j$ then for all $k \in \mathrm{int}(\pi) \neq \varnothing$ the graph $K_n[\pi \setminus \{k\}]$ contains the odd $(i,j)$-path $\pi'=(i,j)$, hence $\pi$ is not minimal in this case. The proof for an even path is similar and omitted.
\qed
\end{Lem}

Given a permutation $\sigma \in S_n$ of $[n]$ and a set $L \subseteq [n]$ let $\succ$ denote the lexicographic monomial order on $S$ induced by
\begin{equation*}
t_{\sigma(1)} \succ  \cdots \succ t_{\sigma(n)} \succ t'_{\sigma(1)} \succ \cdots \succ t'_{\sigma(n)}
\end{equation*}
where $t_{\sigma(i)}=x_{\sigma(i)}$, $t'_{\sigma(i)}=y_{\sigma(i)}$ for all $i \in [n] \setminus L$ and $t_{\sigma(i)}=y_{\sigma(i)}$, $t'_{\sigma(i)}=x_{\sigma(i)}$ for all $i \in L$. For $i,j \in [n]$ write $i \succ j$ if $\sigma^{-1}(i) < \sigma^{-1}(j)$. Let $\gbp{G}$ denote the reduced Gr\"obner basis of $\pbei$ with respect to $\succ$. For a nonzero $f \in S$ let $N_{\mathcal{G}_{\succ}(G)}(f)$ denote the normal form of $f$ with respect to $\gbp{G}$ and let $in_{\succ}(f)$ denote the initial monomial of $f$ with respect to $\succ$. For the next three lemmas (\ref{lemma:reduction1} to \ref{lemma:basis}) fix a permutation $\sigma \in S_n$ and a set $L \subseteq [n]$. For $v \in [n]$
\begin{equation*}
c_v:=
\begin{cases}
+1 &\mbox{if } \sigma^{-1}(v) \notin L \\
-1 &\mbox{if } \sigma^{-1}(v) \in L;
\hspace{5mm}
\end{cases}
r_v:=
\begin{cases}
y_v &\mbox{if } \sigma^{-1}(v) \notin L \\
x_v &\mbox{if } \sigma^{-1}(v) \in L.
\end{cases}
\end{equation*}
For $i,j,k,l \in [n]$ let $B_{(i,j)} = \{c_i(x_ix_j-y_iy_j):i \succ j\}$, $B_{(i,k,j)} = \{c_i(x_iy_j-y_ix_j)r_k:i,k \succ j\}$ and $B_{(i,k,l,i)} = \{c_i(x_i^2-y_i^2)r_kr_l:k,l \succ i\}$. Notice that for an element $f \in \left(B_{(i,j)} \cup B_{(i,k,j)} \cup B_{(i,k,l,i)}\right)$, the value of $c_i$ ensures that the coefficient of the initial monomial $in_{\succ}(f)$ is $1$. Finally let the set $\Gamma \subseteq B_{(i,k,j)}$ consist of all binomials $f = c_i(x_iy_j-y_ix_j)r_k \in B_{(i,k,j)}$ satisfying $i \succ k \succ j$ and $|\{\sigma^{-1}(i),\sigma^{-1}(k)\} \cap L|=1$.

\begin{Lem}\label{lemma:reduction1}
Let $f \in S$ be a nonzero binomial corresponding to a minimal path $\pi$ in $K_n$ (in the sense of Lemma~\ref{lemma:ksw15}). Then $f$ is reduced with respect to $\succ$ if and only if $f \in \Lambda := \left(B_{(i,j)} \cup B_{(i,k,j)} \cup B_{(i,k,l,i)}\right) \setminus \Gamma$.
\proof By Lemma~\ref{lemma:minimal} it suffices to consider only binomials corresponding to the paths $(i,j)$, $(i,k,j)$ and $(i,k,l,i)$ in $K_n$, where $i,j,k$ and $l$ are distinct elements of $[n]$. Without loss of generality we may assume that $i \succ j$. If $f$ is the binomial corresponding to the path $(i,j)$ then clearly $f$ is reduced if and only if $f \in B_{(i,j)} \subseteq \Lambda$.

If $f$ is a binomial corresponding to the path $(i,k,j)$ then there are three conceivable cases: $i \succ j \succ k$, $i \succ k \succ j$ and $k \succ i \succ j$. For ease of notation $f_{(i_1,i_2,i_3)}^t:=c_{i_1}(x_{i_1}y_{i_3}-y_{i_1}x_{i_3})t_{i_2}$ where the superscript $t \in \{x,y\}$ indicates whether $t_{i_2}=x_{i_2}$ or $t_{i_2}=y_{i_2}$.

\begin{Case}[$i \succ j \succ k$]
If $\sigma^{-1}(i) \not\in L$ then $f_{(i,k,j)}^x$ is reduced by $x_ix_k-y_iy_k$; also $f_{(i,k,j)}^y$ is reduced by $f_{(i,j,k)}^y$ if $\sigma^{-1}(j) \not\in L$, or $y_jy_k-x_jx_k$ if $\sigma^{-1}(j) \in L$. If $\sigma^{-1}(i) \in L$ then $f_{(i,k,j)}^x$ is reduced by $x_jx_k-y_jy_k$ if $\sigma^{-1}(j) \not\in L$, or $f_{(i,j,k)}^x$ if $\sigma^{-1}(j) \in L$; also $f_{(i,k,j)}^y$ is reduced by $y_iy_k-x_ix_k$.
\end{Case}

\begin{Case}[$i \succ k \succ j$]
If $\sigma^{-1}(i) \not\in L$ then $f_{(i,k,j)}^x$ is reduced by $x_ix_k-y_iy_k$; also $f_{(i,k,j)}^y$ is reduced by $y_ky_j-x_kx_j$ if $\sigma^{-1}(k) \in L$ but is irreducible if $\sigma^{-1}(k) \not\in L$. If $\sigma^{-1}(i) \in L$ then $f_{(i,k,j)}^x$ is reduced by $x_kx_j-y_ky_j$ if $\sigma^{-1}(k) \not\in L$ but is irreducible if $\sigma^{-1}(k) \in L$; also $f_{(i,k,j)}^y$ is reduced by $y_iy_k-x_ix_k$.
\end{Case}

\begin{Case}[$k \succ i \succ j$]
If $\sigma^{-1}(i) \not\in L$ then $f_{(i,k,j)}^x$ is reduced by $x_kx_i-y_ky_i$ if $\sigma^{-1}(k) \not\in L$ but is irreducible if $\sigma^{-1}(k) \in L$; also $f_{(i,k,j)}^y$ is reduced by $y_ky_j-x_kx_j$ if $\sigma^{-1}(k) \in L$ but is irreducible if $\sigma^{-1}(k) \not\in L$. If $\sigma^{-1}(i) \in L$ then $f_{(i,k,j)}^x$ is reduced by $x_kx_j-y_ky_j$ if $\sigma^{-1}(k) \not\in L$ but is irreducible if $\sigma^{-1}(k) \in L$; also $f_{(i,k,j)}^y$ is reduced by $y_ky_i-x_kx_i$ if $\sigma^{-1}(k) \in L$ but is irreducible if $\sigma^{-1}(k) \not\in L$.
\end{Case}

Finally let $f=c_i(x_i^2-y_i^2)t_kt_l$ be a binomial corresponding to the path $(i,k,l,i)$. If $i \succ k$ then $f$ is reduced by one of $f_{(i,l,k)}^x, f_{(i,l,k)}^y$ or $c_i(x_ix_k-y_iy_k)$. The case $i \succ l$ is similar. If $k,l \succ i$ then by arguing as before one finds that $f$ is irreducible if and only if $t_s=y_s$ whenever $\sigma^{-1}(s) \not\in L$ and $t_s=x_s$ whenever $\sigma^{-1}(s) \in L$.
\qed
\end{Lem}

\begin{Lem}\label{lemma:reduction2}
Let $\pi$ be an $(i,j)$-path in $K_n$ and $t_k \in \{x_k,y_k\}$ arbitrary. Then $(x_ix_j-y_iy_j)\Pi_{k \in \mathrm{int}(\pi)}t_k$ if $\pi$ is odd and $(x_iy_j-y_ix_j)\Pi_{k \in \mathrm{int}(\pi)}t_k$ if $\pi$ is even, reduce to zero modulo $\Lambda$.
\proof It suffices to restrict to a minimal path $\pi$ (if $\pi$ is not minimal, then its binomial is a multiple of the binomial for a shorter path). If $\pi$ is minimal, then Lemma~\ref{lemma:reduction1} gives the result.
\qed
\end{Lem}

\begin{Lem}\label{lemma:basis}
$\gbp{K_n}=\Lambda$.
\proof The proof is by Buchberger's criterion and is similar to the proof of Theorem 3.6 in Kahle et al.~\cite{ksw:parity_binomial_edge_ideals}. Let $g,g' \in \Lambda$ be reduced binomials corresponding, respectively, to the odd path $\pi=(v_0,v_1,v_2,v_0)$ in $K_n$ and the even path $\pi'=(u_0,u_1,u_2)$ in $K_n$ with $u_0 \succ u_2$. We write
\begin{align*}
		g&=
		c_{v_0}(x_{v_0}^2-y_{v_0}^2)r_{v_1}r_{v_2}, \\
		g'&=
		c_{u_0}(x_{u_0}y_{u_2}-y_{u_0}x_{u_2})r_{u_1}.
\end{align*}
Consider the case that $u_0=v_0$. If $\sigma^{-1}(v_0) \in L$ then 
\begin{equation*}
		\mathrm{spol}(g,g')=
		(y_{v_0}y_{u_2}-x_{v_0}x_{u_2})x_{v_0} 
		\cdot \mathrm{lcm}(r_{v_1}r_{v_2},r_{u_1})
\end{equation*}
which is a monomial multiple of the binomial corresponding to the odd path $(v_0,u_2)$ in $K_n$. Thus $\mathrm{spol}(g,g')$ reduces to zero by Lemma~\ref{lemma:reduction2}. The subcase $\sigma^{-1}(v_0) \not\in L$ is dual to this. In a similar fashion all $\mathrm{spol}(g,g')$ (where $g,g' \in \Lambda$) reduce to zero with respect to $\Lambda$. Thus the set $\Lambda$ fulfills Buchberger's criterion and hence is a Gr\"obner basis of $\mathcal{I}_{K_n}$. By Lemma~\ref{lemma:reduction1} it follows that the elements of $\Lambda$ are reduced with respect to $\succ$.
\qed 
\end{Lem}

\begin{Lem}\label{lemma:normalform0}
Let $f=\mon{u}-\mon{v} \in S$ be multi-homogeneous with $\vertices(\mon{u})=\{i,j\}$. If $\mathrm{gcd}(\mon{u},\mon{v})=1$ then $f \in \{\pm(x_i^px_j^q-y_i^py_j^q),\pm(x_i^py_j^q-y_i^px_j^q):p,q\geq1, \ \mathrm{parity}(p)=\mathrm{parity}(q)\}$.
\proof
By homogeneity $f$ is necessarily of the form
\begin{equation*}
		\mon{u}-\mon{v}=
		x_i^{d_i}y_i^{e_i}x_j^{d_j}y_j^{e_j}-
		x_i^{d_i'}y_i^{e_i'}x_j^{d_j'}y_j^{e_j'}
\end{equation*}
satisfying $d_i+e_i=d_i'+e_i'$ and $d_j+e_j=d_j'+e_j'$. By $\mathrm{gcd}(\mon{u},\mon{v})=1$ it follows that if $d_i>0$ then $d_i'=0$ and similarly for all exponents. But if $d_i>0$ then  $d_i'>0$ or $e_i'>0$ i.e.\ $e_i'>0$ (since $d_i'=0$) which in turn implies $e_i=0$. If $d_j>0$ we get a similar result. If $d_j=0$ then by $\vertices(\mon{u})=\{i,j\}$ we have $e_j>0$. By inverting the argument we obtain
\begin{equation*}
		\mon{u}-\mon{v} \in 
		\{\pm(x_i^{d_i}x_j^{d_j}-y_i^{e_i'}y_j^{e_j'}),
		\pm(x_i^{d_i}y_j^{e_j}-y_i^{e_i'}x_j^{d_j'})\}.
\end{equation*}
The result follows from the implications of homogeneity.
\qed
\end{Lem}

\begin{Lem}\label{lemma:normalform1}
Let $\succ$ be the lexicographic monomial order on $S$ corresponding to $\sigma = \mathrm{id}$ and $L=\varnothing$. Let $\mon{u} = x_1^{d_1}y_1^{e_1} \cdots x_n^{d_n}y_n^{e_n} \in S$ where $|\vertices(\mon{u})|>2$. Let $k=\mathrm{max}\{i:i \in \vertices(\mon{u})\}$ and let $\gamma = \sum_{i=1}^k d_i$. 
\begin{equation}\label{eq:normalform1}
N_{\gbp{K_n}}(\mon{u}) = 
\begin{cases}
y_1^{d_1+e_1} \cdots y_k^{d_k+e_k}, & \mathrm{if} \ \gamma \ \mathrm{is \ even}\\
x_ky_1^{d_1+e_1} \cdots y_{k-1}^{d_{k-1}+e_{k-1}}y_k^{d_k+e_k-1}, & \mathrm{if} \ \gamma \ \mathrm{is \ odd}.
\end{cases}
\end{equation}
\proof 
The forms in \eqref{eq:normalform1} are clearly irreducible. By Lemma~\ref{lemma:basis} we have $x_ix_j-y_iy_j \in \gbp{K_n}$ for all $i \succ j\in [n]$, so that $\mon{u}$ can be reduced to
\begin{equation*}
		\mon{u'} =
		x_s^{d_s-l}y_1^{d_1+e_1} \cdots 
		y_{s-1}^{d_{s-1}+e_{s-1}}y_s^{e_s+l}y_{s+1}^{d_{s+1}+e_{s+1}} 
		\cdots y_k^{d_k+e_k}
\end{equation*}
for some $1 \leq s \leq k$ where $l \in \mathbb{Z}_{\geq 0}$, $l\leq d_s$. By the homogeneity of $x_ix_j-y_iy_j$ parity($\sum_{i=1}^k d_i$) = parity($d_s-l$). If $d_s=l$ then we are done. Otherwise we consider the following two cases. If $1 \leq s \leq k-1$ then by $|\vertices(\mon{u})|>2$ there exists $v \in \vertices(\mon{u}) \setminus \{s\}, \ v \neq k$. Since $v,s \succ k$, by Lemma~\ref{lemma:basis} $f = y_v(x_sy_k-y_sx_k) \in \gbp{K_n}$. Using $f$ and $x_sx_k-y_sy_k \in \gbp{K_n}$ in that order $\mon{u'}$ can be reduced to
\begin{equation*}
		\mon{u''} = 
		x_s^{d_s-l-2}y_1^{d_1+e_1} \cdots 	y_{s-1}^{d_{s-1}+e_{s-1}}y_s^{e_s+l+2}y_{s+1}^{d_{s+1}+e_{s+1}} 
		\cdots y_k^{d_k+e_k}.
\end{equation*}
Repeated iteration of this step gives one of the forms in \eqref{eq:normalform1}, depending on the parity of $d_s-l$. If $s=k$ then by $|\vertices(\mon{u})|>2$ there exist $v_1,v_2 \in \vertices(\mon{u}) \setminus \{k\}, \ v_1 \neq v_2$. Since $v_1,v_2 \succ k$, by Lemma~\ref{lemma:basis} $f = y_{v_1}y_{v_2}(x_k^2-y_k^2) \in \gbp{K_n}$. Using $f$ we can reduce $\mon{u'}$ to one of the forms in \eqref{eq:normalform1}, depending on the parity of $d_s-l$.
\qed
\end{Lem}

\begin{Lem}\label{lemma:normalform2}
Let $\succ$ be the lexicographic monomial order on $S$ corresponding to $\sigma = \mathrm{id}$ and $L=\varnothing$. Let $\mon{u} = x_i^{d_i}y_i^{e_i}x_j^{d_j}y_j^{e_j} \in S$ i.e. $|\vertices(\mon{u})| \leq 2$. Let $q=\mathrm{min}\{d_i,d_j\}$.
\begin{equation}\label{eq:normalform2}
		N_{\gbp{K_n}}(\mon{u}) = 
		x_i^{d_i-q}y_i^{e_i+q}x_j^{d_j-q}y_j^{e_j+q}.
\end{equation}
\proof 
The monomial $\mon{u}$ can be reduced to \eqref{eq:normalform2} by the binomial $x_ix_j-y_iy_j \in \gbp{K_n}$ (Lemma~\ref{lemma:basis}). No element of the set $\{in_{\succ}(g) : g \in \gbp{K_n}\}$ divides \eqref{eq:normalform2} hence \eqref{eq:normalform2} is irreducible.
\qed
\end{Lem}

\begin{Lem}\label{lemma:normalform3}
Let $f=\mon{u}-\mon{v} \in S$ be such that $|\vertices(\mon{u})|>2$. Then $f \in \mathcal{I}_{K_n}$ if and only if $f$ is multi-homogeneous.
\proof
``$\Rightarrow$'': This is clear. ``$\Leftarrow$'': Suppose that $f$ is multi-homogeneous in degree $((\alpha_1,\alpha_2),\beta) \in \ints^2_2 \times \nats^n$. Let $\succ$ be the lexicographic monomial order on $S$ corresponding to $\sigma = \mathrm{id}$ and $L=\varnothing$. Let $k=\mathrm{max}\{i:i \in \vertices(\mon{u})\}$. If $\alpha_1=0$ then $N_{\gbp{K_n}}(\mon{u})=N_{\gbp{K_n}}(\mon{v})$ is the first form in \eqref{eq:normalform1} and if $\alpha_1=1$ then $N_{\gbp{K_n}}(\mon{u})=N_{\gbp{K_n}}(\mon{v})$ is the second form in \eqref{eq:normalform1}.
\qed
\end{Lem}

\begin{Lem}\label{lemma:normalform4}
Let $f=\mon{u}-\mon{v} \in S$ be such that $|\vertices(\mon{u})| \leq 2$. Then $f \in \mathcal{I}_{K_n}$ if and only if $f$ is of the form $x_i^{d_i}y_i^{e_i}x_j^{d_j}y_j^{e_j}-x_i^{d_i-q}y_i^{e_i+q}x_j^{d_j-q}y_j^{e_j+q}$ where $q \in \ints$.
\proof
``$\Rightarrow$'': $f=\mon{u}-\mon{v} \in \mathcal{I} \Rightarrow N_{\gbp{K_n}}(\mon{u})=N_{\gbp{K_n}}(\mon{v})$ where $\succ$ is the lexicographic monomial order on $S$ corresponding to $\sigma = \mathrm{id}$ and $L=\varnothing$. Apply Lemma~\ref{lemma:normalform2}. ``$\Leftarrow$'': Apply Lemma~\ref{lemma:normalform2}.
\qed
\end{Lem}

\begin{Thm}\label{theorem:parity}
Conjecture~\ref{conjecture:one} holds for $G=K_n$. 
\proof
Let $G=K_n$ throughout. We prove the theorem in three steps; the containments $\sset{\pbei} \subseteq \ugb{\pbei}$, $\ugb{\pbei} \subseteq \gra{\pbei}$ and $\gra{\pbei} \subseteq \sset{\pbei}$.

\setcounter{Step}{0}
\begin{Step}
$\sset{\pbei} \subseteq \ugb{\pbei}$: Here we invoke Lemmas~\ref{lemma:minimal} and \ref{lemma:basis}. If $\pi$ is a path of the form $(i,k,l,i)$ in $G$ then let $\sigma \in S_n$ be a permutation such that $k,l \succ i$. A suitable choice of $L \subseteq \{\sigma^{-1}(i),\sigma^{-1}(k),\sigma^{-1}(l)\}$ then provides that $(x_i^2-y_i^2)r_kr_l \in \gbp{G}$ or $(y_i^2-x_i^2)r_kr_l \in \gbp{G}$. The cases $\pi=(i,k,j)$ and $\pi=(i,j)$ are similar and omitted.
\end{Step}

\begin{Step}
$\ugb{\pbei} \subseteq \gra{\pbei}$: \cite[Proposition 4.2]{ctv:binomial_fibers_indispensable_binomials}.
\end{Step}

\begin{Step}
$\gra{\pbei} \subseteq \sset{\pbei}$: Let $f=\mon{u}-\mon{v} \in \gra{\pbei}$. If $|\vertices(\mon{u})| \leq 2$ then by Lemma~\ref{lemma:normalform4} we have $f = x_i^{d_i}y_i^{e_i}x_j^{d_j}y_j^{e_j}-x_i^{d_i-q}y_i^{e_i+q}x_j^{d_j-q}y_j^{e_j+q}$ where $q \in \ints$. If $q>0$ then necessarily $f=x_ix_j-y_iy_j \in \sset{\pbei}$. If $q<0$ then necessarily $f=y_iy_j-x_ix_j \in \sset{\pbei}$.

Now consider the case that $|\vertices(\mon{u})|>2$. First suppose that $f=t_k(\mon{u'}-\mon{v'})$ where $k \in \vertices(\mon{u})$ and $t_k \in \{x_k,y_k\}$. Now $\mon{u'}-\mon{v'}$ is multi-homogeneous and it must be that $|\vertices(\mon{u'})|=2$ since otherwise by Lemma~\ref{lemma:normalform3} $\mon{u'}-\mon{v'} \in \pbei$, contradicting the primitivity of $f$. Write
\begin{equation}\label{eq:last}
		f=t_k(\mon{u'}-\mon{v'})=
		t_k(x_i^{d_i}y_i^{e_i}x_j^{d_j}y_j^{e_j}-
		x_i^{d_i'}y_i^{e_i'}x_j^{d_j'}y_j^{e_j'}).
\end{equation}
If $\mathrm{gcd}(\mon{u'},\mon{v'}) \neq 1$ then it follows from the previous argument that $f=t_kt_l(x_s^{d_s}y_s^{e_s}-x_s^{d_s'}y_s^{e_s'})$ where $\{s,l\} \subseteq \{i,j\}$ and $t_l \in \{x_l,y_l\}$. Now $x_s^{d_s}y_s^{e_s}-x_s^{d_s'}y_s^{e_s'}$ is multi-homogeneous and nonzero. These criteria are minimally satisfied by $(d_s,e_s) \in \{(2,0),(0,2)\}$ i.e.\ $f=t_kt_l(x_s^2-y_s^2) \in \sset{\pbei}$ or $f=t_kt_l(y_s^2-x_s^2) \in \sset{\pbei}$. For $d_s+e_s>2$ the primitivity of $f$ is contradicted either by one of these binomials or by an element of the form $\pm(x_ix_j-y_iy_j) \in \sset{\pbei}$. If instead in \eqref{eq:last} $\mathrm{gcd}(\mon{u'},\mon{v'})=1$ then by Lemma~\ref{lemma:normalform0} and the primitivity of $f$ we have $f=\pm t_k(x_iy_j-y_ix_j) \in \sset{\pbei}$.

Suppose now that $f=\mon{u}-\mon{v}$ cannot be written as $t_k(\mon{u'}-\mon{v'})$ where $k \in \vertices(\mon{u})$ and $t_k \in \{x_k,y_k\}$. Since $f$ is multi-homogeneous and $|\vertices(\mon{u})|>2$ we can assume that for some $i,j \in \vertices(\mon{u})$ either $x_ix_j|\mon{u}$ and $y_iy_j|\mon{v}$ or $y_iy_j|\mon{u}$ and $x_ix_j|\mon{v}$ i.e. in this case the primitivity of $f$ is contradicted by an element of the form $\pm(x_ix_j-y_iy_j) \in \sset{\pbei}$.
\qed
\end{Step}
\end{Thm}

\noindent {\bf Acknowledgements.} The authors would like to thank the anonymous referee for their careful reading and helpful suggestions.

\bibliographystyle{amsalpha}
\bibliography{bibfile}

\end{document}